\documentclass[11pt]{amsart}
\usepackage{amssymb,amscd}

\newtheorem{thm}{Theorem}[section]

\newtheorem{cor}[thm]{Corollary}

\newtheorem{prop}[thm]{Proposition}
\newtheorem{lemma}[thm]{Lemma}

\theoremstyle{remark}

\theoremstyle{definition}
\newtheorem{defn}[thm]{Definition}

\numberwithin{equation}{section}

\renewcommand{\bar}{\overline}

\newcommand{\F}{{\mathbb{F}}}
\newcommand{\Q}{{\mathbb{Q}}}
\newcommand{\Z}{{\mathbb{Z}}}

\newcommand{\Aut}{\mathrm{Aut}}

\newcommand{\Gal}{\mathrm{Gal}}
\newcommand{\GL}{\mathrm{GL}}
\newcommand{\Frob}{\mathrm{Frob}}
\newcommand{\Span}{\mathrm{Span}}

\newcommand{\tor}{\mathrm{tor}}
\newcommand{\tf}{\mathrm{tf}}
\newcommand{\ab}{\mathrm{ab}}

\newcommand{\im}{\mathrm{im}\;}
\newcommand{\gp}{mathfrak p}
\newcommand{\TK}{K^T}

\begin{document}
\title{A Mordell-Weil theorem for abelian varieties over fields generated by torsion points}
\author{Michael Larsen}
\email{larsen@math.indiana.edu}
\address{Department of Mathematics\\
    Indiana University \\
    Bloomington, IN 47405\\
    U.S.A.}

\thanks{The author was partially supported by NSF grants  DMS-0100537 and DMS-0354772.}

\begin{abstract}
Let $A$ be an abelian variety over a number field, $T_\ell$ the $\ell$-adic Tate module, and
$G_\ell$ the image of the Galois action on $T_\ell$.  Then $H^i(G_\ell,T_\ell)$ is a finite $\ell$-group
which vanishes for $\ell\gg 0$.  We apply this bound for $i=1$ and $i=2$ to show that if $K^{\tor}$ denotes the field generated by all torsion points of $A$, then $A(K^{\tor})$ is the direct sum of its torsion group and a free abelian group.
\end{abstract}
\maketitle

\section{Introduction}
\label{s:intro}

Let $A$ be an abelian variety over a number field $K$ and $F$ an algebraic extension of $K$.   The quotient $A(F)/A(F)_{\tor}$ is always torsion-free and abelian.  When $F$ is a number field, the Mordell-Weil theorem asserts that it is finitely generated 
and therefore free.  At the opposite extreme, when $F$ is algebraically closed, the quotient is divisible and is therefore a $\Q$-vector space.  In this paper, we show that when $F$ is any field generated over 
$K$ by the coordinates of torsion points on $A$, the quotient is again a free abelian group, though it may not be finitely generated.  The idea is that dividing points of infinite order in $A(K)$ gives rise to fields which are in some weak sense linearly disjoint from the fields obtained by dividing torsion points.
This is essentially a question of Kummer theory of the abelian variety and is closely related to ideas of Ba\v smikov \cite{Bas}, Ribet \cite{Ribet}, and Jacquinot-Ribet \cite{JR}.  

Let $T_\ell$ denote the $\ell$-adic Tate module of $A$ over a number field $K$.  
Let $G_\ell$ denote the image of the Galois group in $\Aut(T_\ell)$.   
The crucial point in this paper is that for all $i\ge 0$, $H^i(G_\ell,T_\ell)$ is finite,
and when $\ell\gg 0$, it vanishes for all $i$.  In particular, this is true for $i=1$ and $i=2$.   

Throughout the paper, $K$ is a number field, $\bar K$ an algebraic closure, and $G_K = \Gal(\bar K/K)$.
We write $\rho_\ell$ for any continuous $\ell$-adic representation of $G_K$.  
(In the later part of the paper, $\rho_\ell$ comes from the action of $G_K$ on the $\ell$-adic Tate module $T_\ell$ of an abelian variety $A/K$.)  The image of $\rho_\ell$ is always denoted $G_\ell$.  It is a closed subgroup of $\GL_n(\Q_\ell)$ for some $n$ and is therefore an $\ell$-adic Lie group.  The fixed field of the kernel of $\rho_\ell$ is denoted $K^\ell$.  Group cohomology is always understood to mean continuous cohomology.

\section{Galois cohomology}
\label{s:galcoh}

Let $G$ be a compact $\ell$-adic Lie group and $L$ the Lie algebra of $G$ (\cite{Lazard}~V~2.4.2).  
Let $V$ be a finite-dimensional vector space over $\Q_\ell$ on which $G$ operates
continuously.   There is a corresponding ``infinitesimal'' action of $L$ on $V$ (\cite{Lazard}~V~2.4.6).  By a theorem of D.~Lazard (\cite{Lazard}~V~2.4.10~(ii)),
if $H$ is a small enough open subgroup of $G$,
\begin{equation}
\label{e:small-enough}
H^p(H,V)\tilde\to H^p(L,V),
\end{equation}
where the right hand side denotes Lie algebra cohomology.

\begin{prop}
\label{p:Lie-vanishing}
Let $G$ be a compact topological group and $V$ a finite dimensional vector space over $\Q_\ell$
on which $G$ admits a faithful and semisimple action.  Suppose that no subrepresentation
of $V$ factors through a finite quotient of $G$.  Then for all $i\ge 0$,
$$H^i(G,V) = 0.$$
\end{prop}

\begin{proof} 
Let $H$ be an open normal subgroup of $G$.  The restriction of a semisimple representation to such a subgroup is again semisimple, and
by hypothesis, $V^H = 0$.
Suppose that for some integer $i$ and some normal open subgroup $H\triangleleft G$,
we have $H^q(H,V) = 0$ for all $q\le i$.
The Lyndon-Hochschild-Serre spectral sequence for modules finite dimensional over $\Q_\ell$
(which follows from the usual spectral sequence with finite coefficients together with \cite{Tate}~Corollary 2.2 and \cite{Tate}~Proposition 2.3) asserts
$$E_2^{pq} = H^p(G/H,H^q(H,V))\Rightarrow H^{p+q}(G,V).$$
Thus, $H^i(G,V) = 0$ for this value of $i$.  To prove $H^i(G,V)= 0$ for all $i\ge 0$, we may
therefore assume $G=H$ is such that (\ref{e:small-enough}) holds for all $p$.  Every 
$G$-subrepresentation is also an $L$-subrepresentation, and every $L$-subrepresentation comes from a $H$-subrepresentation for some open subgroup $H$ of $G$ (\cite{Lazard}~V~2.4.6).   Without loss of generality, therefore, we may assume that the $G$-subrepresentations and $L$-subrepresentations of $V$ are in one-to-one correspondence and therefore $V$ is a semisimple representation of $L$.  Applying (\ref{e:small-enough}) for $p=0$, we see that
$V^H = 0$ implies $V^L = 0$.  The proposition now follows from the vanishing of the cohomology of a reductive Lie algebra $L$ with coefficients in a semisimple module $V$ such that
$V^L = 0$ \cite{HS}~Th.~10.
\end{proof}

\begin{prop}
\label{t:rational-vanishing}
Let $K$ be any number field,  $\ell$ a rational prime, $V$ a finite dimensional $\Q_\ell$-vector space, and
$$\rho_\ell\colon G_K\to\Aut(V)$$
a continuous semisimple representation which is pure of weight $\neq 0$.   Let 
$G_\ell = \im\rho_\ell$.
Then for all $i\ge 0$,
\begin{equation}
\label{e:Ql-coeffs}
H^i(G_\ell,V) = 0
\end{equation}
\end{prop}

\begin{proof}
We claim that $V^H\neq 0$ for all open subgroups $H$ of $G_\ell$.  
Indeed, let $\gp$ be any prime ideal of $K$ which is unramified in $\rho_\ell$
and whose Frobenius lies in $H$.
As $\rho_\ell$ is pure of non-zero weight, $\rho_\ell(\Frob_\gp)$ cannot have 1 as eigenvalue.
Therefore, $V^H = 0$.
Applying Proposition~\ref{p:Lie-vanishing} with $G=G_\ell$, we obtain 
(\ref{e:Ql-coeffs}).
\end{proof}

\begin{prop}
\label{p:l-adic-finite}
Let $T$ be a free $\Z_\ell$-module on which $G_K$ acts continuously in such a way that $V = T\otimes\Q_\ell$ satisfies the hypotheses of Theorem~\ref{t:rational-vanishing}.  Then 
$H^i(G_\ell,T)$ is a finite $\ell$-group for all $i\ge 0$.
\end{prop}

\begin{proof}
By Lyndon-Hochschild-Serre, it suffices to check that $G_\ell$ has an open subgroup
$H$ such that $H^q(H,T)$ is finite for all $q\le i$.   
By \cite{Lazard}~III~3.1.3, $G_\ell$ contains an open subgroup $H$ which is
equi-$\ell$-valued, and by \cite{Lazard}~V~2.5.8, such a group is a Poincar\'e group in the
sense of \cite{CG}~I~4.5.
In particular, $H^q(H,T/\ell^k T)$ is finite for all $q,k\ge 0$.  From the cohomology sequence
for 
\begin{equation*}
\begin{CD}
0 @>>> T @>\ell>> T@>>> T\otimes\F_\ell @>>> 0
\end{CD}
\end{equation*}
it follows that that $H^q(H,T)\otimes \F_\ell$ is finite for all $q\ge 0$. 
Since
$$H^q(H,T) = \varprojlim_k H^q(H,T/\ell^k T)$$
is an abelian pro-$\ell$ group whose tensor product with $\F_\ell$ is finite, it must be a 
finitely generated $\Z_\ell$-module.
Finally,
$$H^q(H, T)\otimes_{\Z} \Z[1/\ell] = H^q(H,T\otimes\Z[1/\ell]) = H^q(H,V) = 0.$$
We conclude that  $H^q(H,T)$ is finite.

\end{proof}

If $X$ is a non-singular projective variety over a number field $K$, $i$ is a positive integer,
and $\ell$ is a rational prime,
we expect $V=H^i(X,\Q_\ell)$ to satisfy the hypotheses of Theorem~\ref{t:rational-vanishing}.
This is known to be true in a few important cases, notably when $X$ is an abelian variety and $i=1$
(\cite{Faltings}~Satz~3).
It is this case, or more precisely, its dual, which interests us for the rest of the paper.

\begin{thm}
\label{t:big-l-vanishing}
Let $K$ be a number field and $A$ an abelian variety over $K$.  For each rational prime $\ell$
different from the characteristic of $K$, we define $T_\ell$ to be the $\ell$-adic Tate module of $A$
and $G_\ell$ the image of $G_K$ in $\Aut(T_\ell)$. 
Then for all $\ell\gg 0$, $H^i(G_\ell,T_\ell)=0$ for all $i$.
\end{thm}

\begin{proof}
By a theorem of Bogomolov \cite{Bog}, the group $C_\ell$ of homotheties in $G_\ell$ is
open in $\Z_\ell^*$.  By a refinement of Bogomolov's result due to Serre \cite{Serre} \S2, the index of $C_\ell$ in $\Z_\ell^*$ is bounded independently of $\ell$; in particular, for $\ell\gg 0$, 
it contains some element $z_\ell$ not congruent to $1$ (mod $\ell$).  Thus 
$\lim_{k\to\infty} z_\ell^{\ell^k}$ is a central torsion point generating
a non-trivial finite subgroup $Z_\ell$ of prime-to-$\ell$ order.  We apply the 
Lyndon-Hochschild-Serre spectral sequence for $Z_\ell\triangleleft G_\ell$.  Now $H^q(Z_\ell, T_\ell)$ is a $|Z_\ell|$-torsion group
but $|Z_\ell|$ is invertible in $\Z_\ell$, so $H^q(Z_\ell, T_\ell)=0$ for $q\ge 0$.  
Moreover, $H^0(Z_\ell,T_\ell) = 0$, since $Z_\ell$ has non-trivial scalar action on $T_\ell$.
The theorem follows.
\end{proof}

\begin{cor}
\label{c:uniform}
Let $K$ be a number field and $A$ an abelian variety over $K$.   For every $i\ge 0$, there exists a uniform bound $N$ such that 
$$|H^i(G_\ell,T_\ell/\ell^n T_\ell)| < N$$
for all prime powers $\ell^n$.
\end{cor}

\begin{proof}
The cohomology sequence of the short exact sequence
\begin{equation*}
\begin{CD}
0 @>>> T_\ell @>\ell^n>> T_\ell @>>> T_\ell/\ell^n T_\ell @>>> 0
\end{CD}
\end{equation*}
gives
$$|H^i(G_\ell,T_\ell/\ell^n T_\ell)| \le |H^i(G_\ell,T_\ell)||H^{i+1}(G_\ell,T_\ell)|,$$
whence the bound follows from Theorem~\ref{t:big-l-vanishing}.

\end{proof}

\section{Rational points over the field of division points}
\label{s:mordell-weil}

In this section, we apply the computations of the first section to show that the group of points on an abelian variety over the field $K^{\tor}$ generated by all of its torsion is free modulo torsion.   Given Corollary~\ref{c:uniform}, it is not difficult to control $\ell$-divisibility of points on the abelian variety in the towers of number fields arising from $\ell$-power torsion.    What remains is to control $\ell$-divisibility in the towers arising from prime-to-$\ell$ torsion.  

\begin{thm}
\label{t:free}
Let $A$ be an abelian variety over a number field $K$.
Let 
$$K^{\tor}:=\prod_\ell K^\ell,$$
the field generated over $K$ by the coordinates of all torsion points of $A$.   Then
$A(K^{\tor})\cong M \oplus (\Q/\Z)^{2\dim A}$, where $M$ is a free abelian group.
\end{thm}

For any abelian group $M$, we write $M_{\tf}$ for the (torsion-free) quotient
$M/M_{\tor}$.  If $M\subset N$, then $M_{\tor} = N_{\tor}\cap M$, so the natural
homomorphism $M_{\tf}\to N_{\tf}$ is injective.
We use the following criterion for $M_{\tf}$ to be free:

\begin{lemma}
\label{l:freeness}
Let $M$ be an abelian group.  Then $M_{\tf}$ is a free abelian group if and only if
for every finite-dimensional subspace $V\subset M\otimes\Q=M_{\tf}\otimes\Q$, the group 
$\{m\in M_{\tf}\mid m\otimes 1\in V\}$ is finitely generated.
\end{lemma}

\begin{proof}
If $M_{\tf}$ is free, we choose a $\Z$-basis $S$ and write each element of a $\Q$-basis of $V$ as a
$\Q$-linear combination of elements of $S$.  Thus only finitely many elements of $S$ are involved in $V$, so replacing $M$ by the preimage of $\langle S\rangle$, we may assume $M_{\tf}$ is finitely generated.  In this case, every subgroup of $M_{\tf}$ is again finitely generated.

In the other direction, we define a \emph{partial basis} 
of $M_{\tf}$ to be any subset $S\subset M_{\tf}$ which 
maps to a linearly independent set in $M\otimes\Q$ 
and such that for $m\in M_{\tf}$, $m\otimes 1\in \Span(S\otimes 1)$
implies $m\in\langle S\rangle$.  By Zorn's lemma, there exists a maximal partial basis for $M_{\tf}$, which is then necessarily a basis.
\end{proof}

\begin{defn}
If $\Lambda$ is a subgroup of an abelian group $M$
and $\ell$ is a rational prime, we say that
the \emph{$\ell$-saturation} of $\Lambda$ is the subgroup of $M$ 
consisting of elements $m$ such that
$\ell^k m\in \Lambda$ for some non-negative integer $k$.  We say $\Lambda$ is $\ell$-saturated
(in $M$) if
it equals its own $\ell$-saturation.
\end{defn}

\begin{lemma}
\label{l:saturate}
Let  $K^\ell$ denote the fixed field of $G_\ell$.  Then $A(K)_{\tf}$ is of finite index in its $\ell$-saturation in  $A(K^\ell)_{\tf}$.  If $\ell\gg 0$, $A(K)_{\tf}$ is $\ell$-saturated.
\end{lemma}

\begin{proof}
Let $P\in A(K^\ell)$ represent a class belonging to the $\ell$-saturation of $A(K)_{\tf}$, and let
$k$ be the smallest positive integer such that 
$\ell^k [P]\in A(K)_{\tf}$, i.e., such that $\ell^k P\in A(K)+A(K^\ell)_{\tor}$.   
As $A(K^\ell)_{\tor}$ is $\ell$-divisible, without loss of generality we may assume that
$Q := \ell^k P\in A(K)$.
Consider the cohomology sequence of the short
exact sequence of $G_\ell$-modules
\begin{equation*}
\begin{CD}
0 @>>> A[\ell^k] @>>> A(K^\ell) @>\ell^k>> \ell^k A(K^\ell) @>>> 0.
\end{CD}
\end{equation*}
The image of $Q\in H^0(G_\ell,\ell^k A(K^\ell))$ in 
$H^1(G_\ell,A[\ell^k]) = H^1(G_\ell,T_\ell/\ell^k T_\ell)$ has order $\ell^k$; otherwise
$\ell^{k-1} Q = \ell^k R$ for some $R\in A(K^\ell)$, which means $\ell^{k-1} P - R\in A(K^\ell)_{\tor}$
and therefore $\ell^{k-1}[P] = 0$ in $A(K)_{\tf}$, contrary to the definition of $k$.
The lemma now follows from Corollary~\ref{c:uniform}.
\end{proof}

We have already remarked that the torsion in  $A(K^\ell)$ is $\ell$-divisible.
It follows that if $\Lambda\subset A(K^\ell)_{\tf}$ has inverse image $\tilde\Lambda$ in $A(K^\ell)$,
then the index of $\Lambda$ in its $\ell$-saturation in $A(K^{\tor})_{\tf}$ is the same as the index
of $\tilde\Lambda$ in its $\ell$-saturation in $A(K^{\tor})$.   In particular, $\Lambda$ is $\ell$-saturated if and only if $\tilde\Lambda$ is so.  Note also that by Kummer theory, if $L\supset K^\ell$, 
then every Galois extension of $L$ obtained by adjoining coordinates of points which are obtained from points of $A(K^\ell)$
by $\ell$-division is a pro-$\ell$ extension.

We can now prove the main theorem.

\begin{proof}[Proof of Theorem~\ref{t:free}]
By definition of $K^{\tor}$, we have 
$$A(K^{\tor})_{\tor} = A(\bar K)_{\tor}\cong (\Q/\Z)^{2\dim A}.$$
We will prove that $A(K^{\tor})_{\tf}$ is a
free $\Z$-module, and this will imply that the short exact sequence
$$0\to A(K^{\tor})_{\tor}\to A(K^{\tor})\to A(K^{\tor})_{\tf}\to 0$$
splits, which then implies the theorem.  

Without loss of generality, we may replace $K$ by a larger number field $K'$, since 
$A(K^{\tor})_{\tf} \to A({K'}^{\tor})_{\tf}$ is injective, 
and a subgroup of a free abelian group is again free abelian.  
In order that the extensions $K^\ell$ of $K$ have the good linear disjointness property we will need below, we assume that $K$ is large enough for the purposes of \cite{Serre} \S1.
For any fixed finite dimensional $V\subset A(K^{\tor})\otimes\Q$, we choose $K$ large enough
that $V\subset A(K)\otimes\Q$.  Let $\tilde\Lambda = \{P\in A(K)\mid P\otimes 1\in V\}$.
By Lemma~\ref{l:freeness}, it suffices to prove that $\tilde\Lambda$ is of finite index in its
$\ell$-saturation in $A(K^{\tor})$ for all $\ell$ and $\ell$-saturated for all $\ell\gg 0$.

We now fix a prime $\ell$.
By Lemma~\ref{l:saturate},  there exists a finite extension $K_\ell/K$ such that
the $\ell$-saturation of $\tilde\Lambda$ in $A(K^\ell)$ is contained in $A(K_\ell)$; moreover, 
$K_\ell = K$ if $\ell\gg 0$.  
Replacing $K$ by $K_\ell$ and $\tilde\Lambda$ by its $\ell$-saturation in $A(K^\ell)$
before examining $\ell$-saturations in $K^{\tor}$, we may therefore assume without
loss of generality that $\tilde\Lambda$ is already $\ell$-saturated in $A(K^\ell)$.

For every finite set $T$ of primes, let $\TK$ denote the compositum of $K^{\ell'}$ for $\ell'\not\in T$.
By \cite{Serre}~Lemma~1.2.1, there exists a finite set $T$ of primes such that
$$\Gal(\TK/K) = \prod_{\ell'\not\in T} G_{\ell'}.$$
By Serre's uniform version of Bogomolov's theorem \cite{Serre}, if $\ell\gg 0$, 
$G_\ell$ contains a homothety in $\Z_\ell^*$ not congruent to $1$ (mod $\ell$).  
Thus $\Gal(\TK/K)$ has an element $\tau$ mapping to such a homothety but
acting trivially on $\ell'$-torsion for all $\ell'\notin\{\ell\}\cup T$.
We also 
assume $\ell$ is larger than the order of $\prod_{\ell'\in T}\GL_{2\dim A}(\F_{\ell'})$; in particular,
$\ell\not\in T$.  We claim first that $\tilde\Lambda$ is $\ell$-saturated in $\TK$.  If it were not,
we could choose 
$$\sigma\in \Gal(\TK/K^\ell) = \prod_{\ell'\not\in T\cup\{\ell\}} G_{\ell'}$$
and $P\in A(\TK)$ such that $\ell P\in A(K^\ell)$ but $\sigma(P)\neq P$.
As $\sigma$ and $\tau$ act (respectively) trivially and without fixed points on $\ell$-torsion, we have
\begin{equation}
\label{e:commute}
\begin{split}
\tau\sigma^{-1}\tau^{-1}(P)-P 
= \tau(\sigma^{-1}\tau^{-1}(P)) - P \\
= \tau(\sigma^{-1}(\tau^{-1}(P)-P)+\sigma^{-1}(P))-P \\
=\tau(\tau^{-1}(P)-P+\sigma^{-1}(P))-P \\
= \tau(\sigma^{-1}(P)-P)\neq\sigma^{-1}(P)-P,
\end{split}
\end{equation}
so $\sigma\tau\sigma^{-1}\tau^{-1}$ acts non-trivially on $P$.  This is impossible since
$\sigma$ and $\tau$ commute, being supported on complementary index sets in the product
$\prod_{\ell'\in T} G_{\ell'}$.   Now, $K^{\tor}$ is obtained from $\TK$ by adjoining
coordinates of points of order $\prod_{\ell'\in T}\ell'$ in $A$ and then by adjoining
higher $\ell'$-power points for $\ell'\in T$.   A pro-$\ell'$ extension of a field never contains a
non-trivial $\ell$-extension of that field for $\ell'\neq \ell$, and a Galois extension of degree less than $\ell$ cannot contain a non-trivial $\ell$-extension, so the saturation of $\tilde\Lambda$
does not change in passing from $\TK\supset K^\ell$ to $K^{\tor}$.
Thus $\tilde\Lambda$ is $\ell$-saturated in $A(K^{\tor})$.

Finally, we need to show that $\tilde\Lambda$ is of finite index in its $\ell$-saturation in
$A(K^{\tor})$ for \emph{every} $\ell$.  It is enough to show that it is of finite index in its
$\ell$-saturation in $A(\TK K^\ell)$ since the passage from $\TK K^\ell$ to $K^{\tor}$ involves a finite extension followed by a prime-to-$\ell$ profinite extension.  The proof proceeds exactly as before 
except that Bogomolov's theorem \cite{Bog} guarantees only the existence of
$\bar\tau_0$ in $G_\ell$ which is a non-trivial scalar.  By the first theorem in \cite{Serre},
$\Gal(\TK K^\ell/K)$ contains an open subgroup of $G_\ell\times \prod_{\ell'\not\in T\cup\{\ell\}} G_{\ell'}$,
so setting $\bar\tau = \bar\tau_0^k$, we can lift to $\tau\in \Gal(\TK K^\ell/K)$
whose image in $G_\ell$ a homothety not congruent to $1$ (mod $\ell^m$) and whose image in
$\Gal(\TK/K)$ is trivial.  If $\sigma\in \Gal(\TK/K)$ acts non-trivially on $P\in A(\TK)$ such that
$\ell^{m-1}P\not\in A(K^\ell)$ and $\ell^mP\in A(K^\ell)$, we can again apply (\ref{e:commute}) to conclude that $\sigma$ and $\tau$ do not commute, which is impossible.  Thus, $\tilde\Lambda$
is of index $\le \ell^{(m-1)\dim V}$ in its $\ell$-saturation in $A(\TK K^\ell)$.  The theorem follows.

\end{proof}

\begin{cor}
Let $A$ and $B$ be abelian varieties over a number field $K$.
Let $K_A^{\tor}$ denote the extension of $K$ generated by the torsion points of $A$.
If $K\subset L\subset K_A^{\tor}$, then $B(L)_{\tf}$ is free abelian.
\end{cor}

\begin{proof}
Applying Theorem~\ref{t:free} to $A\times B$, we see that $(A\times B)(K_{A\times B}^{\tor})_{\tf}$
is free abelian.  Therefore,
$$B(L)_{\tf}\subset B(K_A^{\tor})_{\tf}\subset (A\times B)(K_{A\times B}^{\tor})_{\tf}$$
is also free abelian.
\end{proof}

We remark that if $A$ is an elliptic curve and $K=\Q$, then $A(K^{\tor})_{\tf}$ can be determined up to isomorphism: it is $\Z^\omega$.  In general, the rank of $A(K^{\tor})$ must be countable; the difficulty lies in proving that it is infinite.  For elliptic curves over $\Q$ written in Weierstrass form, it is clear that points with rational $x$-coordinates are defined over $\Q^{\ab}\subset \Q^{\tor}$, and by the method of \cite{Larsen}, it is easy to see that they generate an infinite-dimensional subspace of $A(\bar\Q)\otimes\Q$.   This raises the following question: is the rank of $A(K^{\tor})$ infinite for all 
abelian varieties over number fields?

\end{document}